\documentclass[a4paper,12pt]{article}

\usepackage[centertags]{amsmath}
\usepackage{amsfonts}
\usepackage{amssymb}
\usepackage{amsthm}
\usepackage{newlfont}

\usepackage{url}

\usepackage[english]{babel}
\usepackage{enumerate}
\usepackage[final]{graphics}
\usepackage{setspace}

\theoremstyle{remark}

\theoremstyle{definition}

\theoremstyle{definition}

     \newcommand {\beq}  {\begin{equation}}
      \newcommand {\eeq}  {\end{equation}}

\author{Alena Aleksenko$^*$, Evgeny Lakshtanov\thanks{Department of
Mathematics, Aveiro University, Aveiro 3810, Portugal.  This work
was supported by \textit{Center for Research
and Development in Mathematics and Applications}(CIDMA) from the ``{\it Funda\c{c}\~{a}o para a
Ci\^{e}ncia e a Tecnologia}'' (FCT), cofinanced by the European
Community Fund FEDER/POCTI} \thanks{e-mail:
lakshtanov@rambler.ru}}

\title{EXAMPLES OF ADMISSIBLE SIMPLIFICATION OF MATHEMATICAL THEORIES}

\begin{document}
\date{}
\maketitle \textbf{Keywords}: Cellular Automata, Turing Machine,
Minesweeper, simplification of mathematical models, questionnaire

"Mathematicians, like physicists, are pushed by a strong
fascination. Research in mathematics is hard, it is intellectually
painful even if it is rewarding, and you wouldn’t do it without
some strong urge. " [D. Ruelle]. We shall give some examples from
our experience, when we were able to simplify some serious
mathematical models to make them understandable by children,
preserving both aesthetic and intellectual value. The latter is in
particularly measured by whether a given simplification allows
setting a sufficient list of problems feasible for school
students.

For a more evident demonstration of our method we chose primary school
students (6-9 year) as a target group. We give examples of Turing
machine, Cellular automata, Minesweeper, Graph theory.

Nowadays, it is thought of importance to introduce current scientific
achievements to children by telling them of the black holes or
demonstrating some impressive chemical experiments. Can this approach
however satisfy us? A. Zvonkin in his acclaimed work [J.Math. Behaviour,
1992] on early child development translates Poincare’s theory on the
role of subconscious into a practical recommendation: {\it questions are more
important than answers.} He depicted an experiment that targeted to find
out whether it is using attractive materials that engaged children in
his lessons or the lessons per se. « Then I say, «All right, I have to
finish the lesson, but you may play with mosaic.» My words are met with
an unanimous yell of indignation, … «No-o, we want a problem!». That’s
how I understood what the truth was. “{\it Children need
intellectual/aesthetical pleasure of full value}. If one of the halves is
absent, the full value is lost, together with the festive feeling. “”
Thus, if you agree with Poincare’s theory on the role of subconscious
work [“Science and Method”, 1908], you would agree that preparing a
model accompanied by a questionnaire for children will further help them
to do research on a higher qualitative level.

The main motivation of the present article is to promote a thesis
that the mathematical scientific community needs for a significant
volume of the internally generated  educational product for
pre-college students. Really, studies that are set up for pupils
by an active researcher, as a rule, contribute to development of
skills essential for successful research activities.
Unfortunately, this aspect is excluded from school curriculum
scope.   Systematical work in this direction is only possible in
case when an actively working mathematician has a part of his
working time (say 30 percents) solely dedicated to this research.

Many of us introduce our  own children to this research, and usually we chose a transphenomenal  Socrates' dialog's style.
It is of  obvious importance both to identify the kernel of these studies  discussed and to  diversify  them so they could be more
easily employed by  the  colleagues   when working with students. Here we present our results describing an attempt to create a way to introduce
our 7-9 years old pupils to several nontrivial mathematical models. Briefly, we have  started
with the  models that originally could not cause a sustained interest of the child. After a  period of time and
series of experiments with pupils we were able to fix a simplified model equipped with a list of questions and exercises.

The most valuable exercises are those  to allow nontrivial routine multiple repetitions.  It is very important for children to gain positive experience of routine exercises. It would allow them to believe/realize that multiple routine exercises in school is a necessary and natural step to the meaningful beauty as in case of a painter's pupil who is ready to grind colours for hours and hours because he is fascinated by the art. Since generation of the quality routine exercises is extremely  important per se, we  formulate a non-formal law:

{\it A problem is solvable by a child if and only if his interest exceeds the amount of effort. }

The exercise can be considered successful even in case when pupils are not able to find a solution themselves.  {\it Children develop when confronted with material that  has already been able to understand, but not yet to reproduce}. Some of my pupils were able to solve all exercises concerning Cellular automata (Model 1 below), however  just an execution of programs for Turing machine (Model 2 below) is done on a peak of their efforts.

\section{Model 1. GAME OF LIFE}

\begin{figure}\label{newtso}
\includegraphics{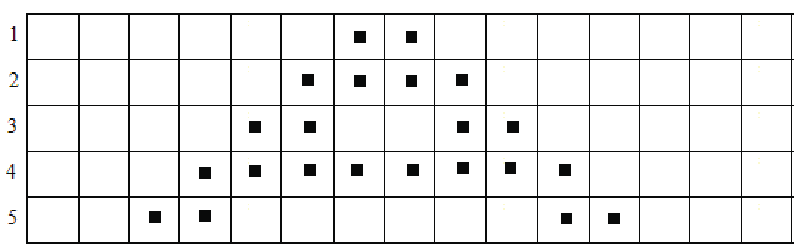} \caption{(a). Example of Dynamics.  }
\end{figure}

Game of Life by John Conway is a well known dynamics with discrete
time on the two-dimensional lattice. Of course, the level of
abstractness is quite high even  to ask children just to find a
descendant. We decided to propose them an excersize with 1D
cellular automata, which inherits the idea of the "Game of Life".
There is a nearest neighbors interaction and life
appears$\backslash$maintained if and only if it has only one
neighbor. The example of the dynamics is presented on the fig.
\ref{newtso}. The great educational advantage with 2D case was
that children could write descendants one below
another.\begin{itemize}
\item Continue a given (local) configuration until the 10th descendant.
\end{itemize}
This exercise is already quite complicated for children, and here we have an effect that Dynamics or Games  have magnetic influent to research it.  When the  pupils reached steady success in solving the exercise one we proposed them to
\begin{itemize}
\item  try to find the ancestor of given local configuration.
\end{itemize}
They did it in the following way:  first, every pupil painted his own local configuration, then he wrote down the descendant on a piece of paper. This piece of paper was transferred to the pupil's neighbor who tried to reconstruct the initial configuration.
\begin{itemize}
\item   Continue a given nonlocal periodical configuration until the 10th descendant.
\item   Find all ancestors of the empty configuration
\item   Find the stationary point of the dynamics
\end{itemize}
Note that in this case  children perform these exercises with
quite abstract objects like {\it inverse map, infinite
configuration, stationary point}. We believe that exercises with
infinite periodical configurations are quite valuable. Here
children have to deal with rather abstract objects. In order to
figure out a descendant they should imagine the cells that are
located beyond their paper sheet. At the same time, such exercises
are both constructive and doable

\begin{figure}\label{life2}
\includegraphics{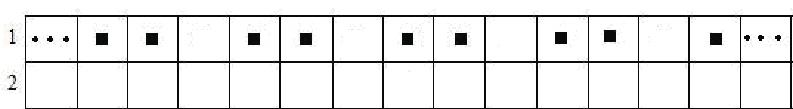} \caption{Stationary configuration. Dots mean periodicity. }
\end{figure}

\section{Model 2. "TURING MACHINE"}

\begin{figure}\label{turing}
\includegraphics{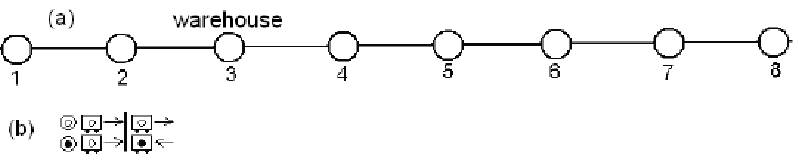} \caption{(a) Tape. (b) Example of set of instructions.}
\end{figure}

\vspace{3mm} For a long time we have been trying to find an
admissible realization of the Turing machine. Once we have read a
juvenile book "EIGHT CHILDREN AND A TRUCK" by Anne-Cath. Vestly.
That is a story about a low-income  family consisted of father,
mother, 8 children and ... a truck. Truck was a quasi alive object
in this story. All of them lived in concord,  and Father and the
truck had to work a lot, as they worked in logistics. I have
decided that this Vestly's world is a good basis for realization
of the Turing machine. "...Previously, they lived in the city, in
a small apartment, but recently they bought a cabin  in the woods,
and could even afford their own cow. A father has a lot of worries
at home, the truck seemed to  become an adult so it could well
work a little bit on its own, to free up Father for his family.
Father has realized  it, they are very well aware of each other,
and he began to think how you could he use it, as the  truck was
still not a human, so he could not  just say: "transfer the goods
there, then go back to the warehouse ". He could remember only one
instruction. Father was a very serious person, who brings
everything to the end, no matter what it takes, we know  this is
not a joke to bring up once the whole eight children… Thus, he
decided to  enumerate all 8 stations that he and the truck were
serving" (Fig. \ref{turing}.(a)). This is an example of the set of
instructions for the truck that it needs to do at a given station
(Fig \ref{turing}.(b)). For example, in the first line it reads:
if you came to the station, moving to the right (arrow), and you
do not have a box in the trunk (the icon with the body blank
circle), and there is nothing for you at the station, then you
have to (look at the icons after the vertical line) to go blank to
the right. A second line reads that if there is a box for you at
the station, then you should pick it up and go to the previous
station. This is an example of the scheme the Father made for each
station.

\begin{itemize}
\item The box is on the warehouse. Transfer it to the station 4.
\item The same, but finally truck should return to the warehouse.
\item {\it Adding.} There are 4 boxes in the warehouse and there are 5 boxes in the 5th station. Truck should transfer all boxes to the warehouse.
Father got the same task on the next day but with another amount of boxes. Should he change the program for the truck?
\item There is one box at the 2nd station and at the 3d stations. Transfer one box to the 4th station and one box to the 5th station.
\item {\it Bisection.} There are 4 boxes at the warehouse. Transfer a half to the 4th station. How these instructions would work if the amount of boxes at the warehouse was 6? What would be the result in case of an odd number of boxes?
\item {\it Subtraction.}  Farther was informed that during the next week, on a daily basis, he should take away from the warehouse that many boxes so that their number becomes equal to the number of boxes at the fourth station.
\item {\it Comparison.} There are boxes at the warehouse and at the 4th station. Determine which station does contain more boxes.
\item {\it Interconversion to the unary system.} There are 5 boxes at the warehouse. Transfer them one by one to the next 5 stations.
\item {\it Invert.}  Once it turned out that the boxes were loaded at the warehouse in the wrong order. So farther asked to pull all of them out and reload them back in the reverse order, so that the first box that was previously put in the warehouse should now go last.
\end{itemize}

\begin{figure}\label{turing2}
\includegraphics{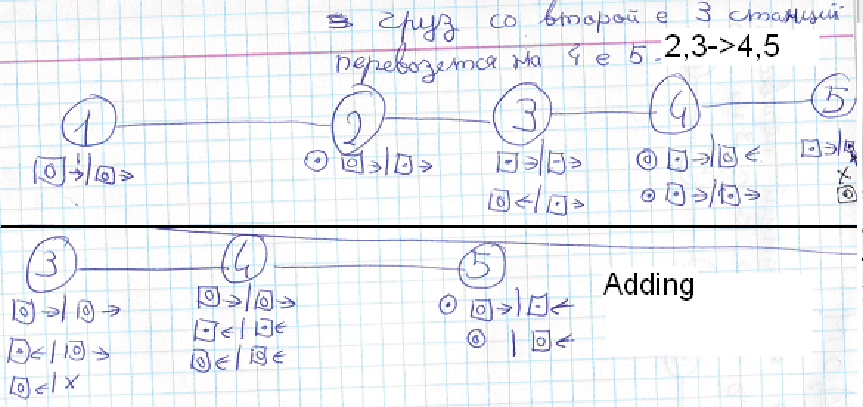} \caption{Example of the programm written by a 9 years old girl.}
\end{figure}

\section{Model 3. MINESWEEPER}

\vspace{3mm}
The following example is not a {\it simplification of the model} in its original form.
"Minesweeper" is a computer game, quite popular among  people that may be very far from science. This was a subject of a number of mathematical researches and not only in terms  of the complexity theory: recently it was proved [R. Kaye] that {\it Turing machine computations can be
reduced to infinite minesweeper problems.}

Minesweeper's player always solves nontrivial but not too
sophisticated tasks. Evidently this game can serve as a training
tool. The problem is  unique, as even starting with  the same
table, players can choose different cells to start with. Thus the
teacher cannot control the difficulty level.  Here it was decided
to develop tables where several cell are open and these open cells
determine the solution uniquely. First tables were produced
manually, however such method proved to be quite time demanding...
Therefore Evgeny Lakshtanov have  started developing an automatic production, and
several months later he has realized that it takes to  analyze
spectral properties of discrete Laplacians in order to prepare a
table for the primary school children to solve.

An algorithm of table generation:
1. Draw a table $m \times n$ similar to a chess board. For every "black" cell Bernoulli test is performed  and if the result is 1, we put a "mine" in this cell. After having run through all the black cells, we fill all "white" cells with numbers (0-4) that are equal to the amount of mines among neighbors.

As a result,  we have a table for the "paper minesweeper" and the distribution of mines constructed in the first step of the algorithm being the solution. How to prove that  this solution is unique?

{\bf Theorem.}\cite[Th. 1]{l}.  If numbers $m+1$ and $n+1$ are coprime, then the constructed table admits one solution only.

\begin{figure}\label{mine3}
\begin{center}
\includegraphics{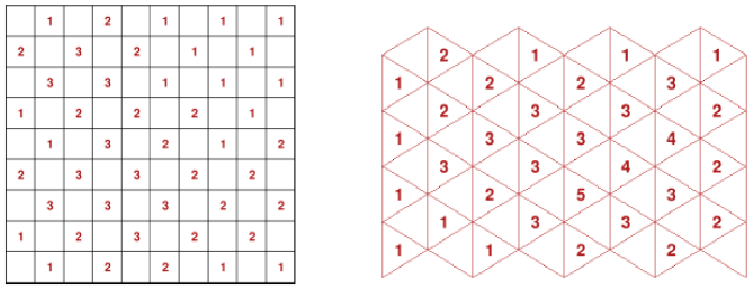} \caption{Example of the table for "paper minesweeper".}
\end{center}
\end{figure}

\section{Own mathematics}
Before these set of exercises we usually ask pupils if they have own experience to introduce objects in mathematics. We talk to them that everyday they study objects like "numbers, functions, triangles" which were introduced much time ago:
(Every pupil need a map of the Europe for next exercises.)

{\bf Definition 1.} Country is {\it monogamous} if it has exactly
one neighbor.
\\
$\bullet$ Find all monogamous country in Europe.
\\
{\bf Definition 2.} A monogamous country is called {\it happy
monogamous} if its neighbor is also monogamous.
\\
$\bullet$ Find all happy monogamous countries in Europe. Why France (Denmark) is not a happy monogamous?
\\
{\bf Definition 3.}  A point is called {\it attractive} if shared
by 3 countries.
\\
$\bullet$ Show several attractive points. (Why this point is not attractive?, Has Denmark any attractive points? Which other countries do not have atractive points?

{\bf Theorem 1.} Monogamous countries do not have attractive
points.
\\
$\bullet$ Prove theorem 1.
\\
$\bullet$. Which attractive point is the closest to Portugal?
\\
$\bullet$. What is the second closest attractive point to
Portugal?

Comment. It is better to use a map where Andorra is colored not
too distinctively, or at least it should be small enough not to
give out the answer too soon that the nearest attractive point to
Portugal also belongs to Andorra.

The next set was particularly chosen to give exercises on the {\it
fit four colours} theorem. Children will need a piece of paper and
a pen/pencil.
\\
{\bf Definition 4.} A country is called {\it friendly} if it has
at least two neighbors, and each two of them also neighbor each
other. Rank of a friendly country is equivalent to a number of its
neighbors.

\begin{figure}\label{map1}
\begin{center}
\includegraphics{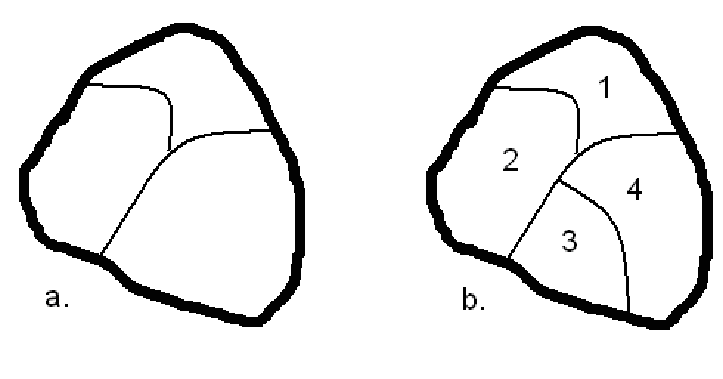} \caption{a. There are 3 friendly countries on the island. }\label{newtsol2}
\end{center}
\end{figure}

$\bullet$ Here is an island in the sea. Check that all 3 countries
on the picture \ref{map1}.a are {\it friendly}.
\\
$\bullet$ What countries on the picture \ref{map1}.b are friendly?
\\
$\bullet$ Paint 4 countries on the island such a way that all of
them should be friendly.

Let us return to the map of Europe.
\\
{Definition 5.} The rank of a friendly country equals the number
of its neighbours.
\\
$\bullet$ Can rank of a friendly country be 1? What is the rank of
Moldova?
\\
$\bullet$ Find all rank 2 friendly countries in Europe. Find all
rank 3 friendly countries in the Europe.
\\
$\bullet +$ Can rank of a friendly country be 4?

Some monkeys use tail as a limb. Now, formally we can think that
they have an odd number of limbs.
\\
{Definition 6.} The country is {tailed} if it has an odd number of
neighbours. A country is {\it tailless} in case this number is
even.
\\
$\bullet +$ Prove that for any map, amount of tailed  countries is
even.

\bigskip

\noindent \textbf{Acknowledgement.} The authors are very grateful
to Oleg German  for numerous and useful
discussions.


\begin{thebibliography}{102}

\bibitem{z} A.K. Zvonkin. {\it Mathematics for little ones.} – Journal of Mathematical
Behavior, 1992, vol. 11, no. 2, 207–219. {\it Children and $C_5^2$.},
1993, vol. 12, no. 2, 141–152.


\bibitem{l} E.Lakshtanov, O.German, {\it Paper Minesweeper" or how to play "Minesweeper" without a computer},  Applicable Analysis, 89(12), 2010, 1907-1916.
\end{thebibliography}
\end{document}